\documentclass[12pt,notitlepage]{article}
\usepackage{amssymb}

\usepackage[dvips]{graphicx}
\usepackage{amsmath}

\newtheorem{theorem}{Theorem}[section]

\newtheorem{corollary}[theorem]{Corollary}

\newtheorem{lemma}[theorem]{Lemma}

\newtheorem{proposition}[theorem]{Proposition}

\newenvironment{proof}[1][Proof]{\textbf{#1.} }{\ \rule{0.5em}{0.5em}}

\begin{document}

\title{On the product of vector spaces in a field extension}
\title{On linear versions of some addition theorems}
\date{}
\author{Shalom Eliahou and C\'{e}dric Lecouvey \\
LMPA Joseph Liouville, FR CNRS 2956\\
Universit\'{e} du Littoral C\^{o}te d'Opale\\
50 rue F. Buisson, B.P. 699\\
F-62228 Calais cedex, France}
\maketitle

\begin{abstract}

Let $K \subset L$ be a field extension. Given $K$-subspaces $A,B$ of $L$, we
study the subspace $\langle AB \rangle$ spanned by the product set $AB=\{ab
\mid a \in A, b \in B\}$. We obtain some lower bounds on $\dim_K \langle AB
\rangle$ and $\dim_K \langle B^n \rangle $ in terms of $\dim_K A$, $\dim_K B$ 
and $n$. This is achieved by establishing linear versions of constructions 
and results in additive number theory mainly due to Kemperman and Olson.
\end{abstract}

\section{Introduction}

Let $G$ be a group, written multiplicatively. Given subsets $A,B \subset G$,
we denote by 
\begin{equation*}
AB = \{ab \mid a \in A, b\in B \}
\end{equation*}
the \textit{product set} of $A,B$. For $A,B$ finite, several results in
additive number theory give estimates on the cardinality of $AB$ as a
function of $|A|,|B|$. Instances of such results are the Cauchy-Davenport
theorem for cyclic groups of prime order, a theorem of Kneser for abelian
groups, and theorems of Kemperman and of Olson for possibly nonabelian
groups.

In this paper, we consider the following analogous question. Given a field
extension $K\subset L$, and finite-dimensional $K$-subspaces $A,B$ of $L$,
define 
\begin{equation*}
\langle AB\rangle =\text{ the $K$-subspace spanned by the product set $AB$
in $L$.}
\end{equation*}
What can be said, then, about the dimension of the subspace $\langle
AB\rangle $? The main object of this paper is to establish linear analogues,
in this new setting, of several results in additive number theory. In
particular, we shall obtain nontrivial lower bounds on $\dim \langle
AB\rangle$ in terms of $\dim A$, $\dim B$.

This question has barely been addressed in the literature. What seems to be
the so far unique result of this type is due to Hou, Leung and Xiang \cite
{Xian} and is given below. We first recall Kneser's theorem from additive
number theory \cite{Kn53, Kn55}.

\begin{theorem}[Kneser]
\label{Kneser} Let $G$ be an abelian group and let $A,B\subset G$ be finite
nonempty subsets. Then 
\begin{equation*}
|AB|\geq |AH|+|BH|-|H|,
\end{equation*}
where $H=\{g\in G\mid gAB=AB\}$ is the stabilizer of $AB$.
\end{theorem}

The following linear version has been obtained in \cite{Xian},
and has motivated us to further explore the links between additive
number theory, linear algebra and field extensions.

\begin{theorem}[Hou, Leung and Xiang]
\label{linear Kneser} \label{TH_HLX} Let $K\subset L$ be a commutative field
extension, and let $A,B\subset L$ be nonzero finite-dimensional $K$-subspaces 
of $L$. Suppose that every algebraic element in $L$ is separable
over $K.$ Let $H$ be the subfield of $L$ which stabilizes $\langle AB\rangle.$ Then 
\begin{equation*}
\dim _{K}\langle AB\rangle \geq \dim _{K}A+\dim _{K}B-\dim _{K}H.
\end{equation*}
\end{theorem}

The paper is organized as follows. Several classical results and
constructions in additive number theory are recalled in Section~\ref{sec-kem}. 
In Section 3, we give a new variant of a theorem of Olson. The switch to
the field extension setting is performed from Section~\ref{sec_linear
setting} on, where we give linear analogues of the additive results stated
in the preceding sections.

\section{Classical addition theorems\label{sec-kem}}

\subsection{A basic result}

One of the simplest results on the product of two sets in a finite group is
the following. Given a subset $X$ of a group $G$, we denote $X^{-1}=\{x^{-1}
\mid x \in X \}$.

\begin{proposition}
\label{basic} Let $G$ be a finite group. Let $A,B$ be nonempty subsets of $G$. 
If $|A|+|B|>|G|$, then $AB=G$.
\end{proposition}

The proof, informally, goes as follows. For $x \in G$, we have $A^{-1}x \cap
B \neq \emptyset$, since $|A^{-1}x|+|B|$ = $|A|+|B| > |G|$. It follows that $x \in AB$.

\smallskip

Linearizing this result requires a little more work. This is done in Section~\ref{lin ols}, with various proofs.

\subsection{Theorems of Kemperman and Olson}

Here are a few classical results in additive number theory, to be linearized in subsequent sections. 
We refer the reader to \cite{ols} for proofs.
\begin{theorem}[Kemperman]
\label{th_ke}\label{th_kem}\label{THK} Let $A,B$ be two finite subsets of a
group $G.\;$Assume there exists an element $c\in AB$ appearing exactly once
as a product $c=ab$ with $a\in A$, $b\in B$. Then 
\begin{equation*}
\left| AB\right| \geq \left| A\right| +\left| B\right| -1.
\end{equation*}
\end{theorem}

The next result is a nonabelian analogue of Kneser's Theorem~\ref{Kneser}.

\begin{theorem}[Olson]
\label{thOl} Let $A,B$ be two finite subsets of a group $G.\;$There exists a
nonempty subset $S$ of $AB$ and a finite subgroup $H$ of $G$ such that 
\begin{equation*}
\left| AB\right| \geq \left| S\right| \geq \left| A\right| +\left| B\right|
-\left| H\right| 
\end{equation*}
and either $HS=S$ or $SH=S$.
\end{theorem}

Here is an immediate corollary.

\begin{theorem}[Kemperman]
\label{torsion-free} Let $G$ be a torsion-free group. Let $A,B$ be 
nonempty finite subsets of $G$. Then 
\begin{equation*}
\left| AB\right| \geq \left| A\right| +\left| B\right| -1.
\end{equation*}
\end{theorem}

\noindent Olson \cite{ols} also derived, as a consequence of Theorem~\ref
{th_kem} and \ref{thOl}, the following two interesting results. 
Linear versions are given in Section~\ref{sec: powers}.

\begin{theorem}[Olson]
\label{thOl2} Let $A,B$ be finite subsets of a group $G$ with $1\in
B$. Then 
\begin{equation*}
AB^{2}=AB\;\text{ or }\;\left| AB\right| \geq \left| A\right| +\frac{1}{2}\left| B\right|.
\end{equation*}
\end{theorem}

\begin{theorem}[Olson]
\label{Th_ol3} Let $B$ be a finite subset of a group $G$. Then 
\begin{equation*}
\left| B^{n}\right| =\left| B^{n+1}\right| \;\text{ or }\;\left|
B^{n}\right| \geq \left| B^{n-1}\right| +\frac{1}{2}\left| B\right| .
\end{equation*}
\end{theorem}

\subsection{A main tool: the Kemperman transform}\label{tool}

The above results are obtained by cleverly iterating Kemperman transforms,
which we now recall.

Let $G$ be a (possibly nonabelian) multiplicative group and $(A,B)$ a pair of nonempty
finite subsets of $G.\;$Let $x$ be any element in $G$.\ The Kemperman
transformed pairs $(A^{\prime },B^{\prime })$ and $(A^{\prime \prime
},B^{\prime \prime })$ with respect to $x$ are defined by 
\begin{equation*}
(A^{\prime}  =  A\cup Ax, \, B^{\prime}  =  B\cap x^{-1}B)
\end{equation*}
and
\begin{equation*}
(A^{\prime\prime}  =  A \cap Ax^{-1},  \, B^{\prime\prime}  =  B\cup xB).
\end{equation*}
The following properties are straightforward to check.

\begin{enumerate}
\item  $A^{\prime }B^{\prime }\subset AB$ \ and\ $A^{\prime \prime
}B^{\prime \prime }\subset AB$;

\item  $B^{\prime }\neq \emptyset $ $\Longleftrightarrow $ $x\in BB^{-1}$\ and\ 
$A^{\prime \prime }\neq \emptyset $ $\Longleftrightarrow $ $x\in A^{-1}A$;

\item  $\left| A^{\prime }\right| >\left| A\right| $ $\Longleftrightarrow $ 
$Ax\neq A$ \ and \ $\left| B^{\prime \prime }\right| >\left| B\right| $ 
$\Longleftrightarrow $ $xB\neq B$.
\end{enumerate}

\begin{proposition}[Kemperman \protect\cite{Kem}]
Let $(A,B)$ be a pair of nonempty finite subsets in a group $G.\;$Assume that there is
an element $d\in A^{-1}A\cap BB^{-1}$ such that either $Ad\neq A$ or $dB\neq
B$. Then there is a pair of nonempty subsets $(A_{1},B_{1})$, obtained as a Kemperman transform of $(A,B)$,
satisfying:

\begin{enumerate}
\item[(1)] $A_{1}B_{1}\subset AB,$

\item[(2)] either $\left| A_{1}\right| +\left| B_{1}\right| =\left| A\right|
+\left| B\right| $ and $\left| A_{1}\right| >\left| A\right| ,$ or $\left|
A_{1}\right| +\left| B_{1}\right| >\left| A\right| +\left| B\right| .$
\end{enumerate}
\end{proposition}

\section{A variant of a theorem of Olson}

\label{variant olson}

In this section we give a variant of Theorem~\ref{thOl2} which is neither stronger
nor weaker, in the sense that neither result implies the other one. 
A linearized version of this variant will be presented in Section~\ref{sec: lin ABC}.

\begin{theorem}
\label{ABC} Let $A,B,C$ be nonempty finite subsets of a group $G$. Suppose 
$B\subset C$ and $1\in C.$ Then 
\begin{equation*}
ABC=AB\;\text{ or }\;\left| ABC\right| \geq \left| A\right| +\left| B\right|
.
\end{equation*}
\end{theorem}

\begin{proof}
The proof is close to that of Theorem \ref{thOl2} in \cite{ols}. If $C=\{1\}$
then $ABC=AB.$ We now assume $|C|\geq 2$ and proceed by induction on $\left|
AB\right| .$

If $\left| AB\right| =1$ then $\left| A\right| =\left| B\right| =1,$ and 
$|ABC|=|C|\geq 2\geq |A|+|B|$.

Assume now $\left| AB\right| >1$ and our statement true for any pair 
$(A_{1},B_{1})$ such that $B_{1}\subset C$ and $\left| A_{1}B_{1}\right|
<\left| AB\right| .$ As in Olson's original proof of Theorem \ref{thOl2}, we
consider two cases.

\medskip 

\noindent $(i)$ Assume there exists a nonempty subset $S\subset AB$ such that 
$SC=S$. Since $S$ is finite, we have $Sc=S$ for all $c\in C$, whence 
$S{C}^{-1}=S.$

Set $A_{0}=A\cap S$ and $A_{1}=A\setminus S,$ so that $A=A_{0}\sqcup A_{1}$.
We claim that 
\begin{equation}
A_{0}B=S.  \label{eq: A0B}
\end{equation}
Indeed, as $B\subset C$, we have $A_{0}B\subset SC=S$. Conversely, let $s\in
S$. Then $s=ab$ for some $a\in A$ and $b\in B$. Since $a=sb^{-1}\in
SC^{-1}=S,$ it follows that $a\in A_{0}$, whence $s\in A_{0}B$ as desired.

As a first consequence, we obtain 
\begin{equation}
AB=S\sqcup A_{1}B.  \label{eq: AB}
\end{equation}
Indeed, we have $AB=A_{0}B\cup A_{1}B=S\cup A_{1}B$. Moreover, as $S\cap
A_{1}=\emptyset $, we have $SB^{-1}\cap A_{1}=\emptyset $ since 
$SB^{-1}\subset SC^{-1}=S$. It follows that $S\cap A_{1}B=\emptyset .$

Our next claim is that 
\begin{equation}
ABC=S\sqcup A_{1}BC.  \label{eq: ABC}
\end{equation}
Indeed, we have $ABC=SC\cup A_{1}BC=S\cup A_{1}BC$. The intersection $S\cap
A_{1}BC$ is indeed empty, since $S=SC^{-1}$ and so $SC^{-1}\cap A_{1}B=S\cap
A_{1}B=\emptyset $.

Consequently, we derive 
\begin{equation}
|ABC|\geq |S|+|A_{1}BC|\geq |A_{0}|+|A_{1}BC|,  \label{eq: cardABC}
\end{equation}
where the estimate $|S|\geq |A_{0}|$ follows from the equality $S=A_{0}B$ in
(\ref{eq: A0B}).

We may assume $A_{1}\not=\emptyset $, for otherwise $S=A_{0}B=AB$, whence 
$ABC=SC=S=AB$ and we are done. We have $|A_{1}B|<|AB|$ by (\ref{eq: AB}). By
our induction hypothesis, we either have $A_{1}BC=A_{1}B$ or $|A_{1}BC|\geq
|A_{1}|+|B|$.

If $A_{1}BC=A_{1}B$, then by (\ref{eq: ABC}) and 
(\ref{eq: AB}), we obtain $ABC=S\cup A_{1}BC=S\cup A_{1}B=AB$ and we are done in this
case.

If now $|A_{1}BC|\geq |A_{1}|+|B|$, then by (\ref{eq: cardABC}) 
and $|A|=|A_{0}|+|A_{1}|$, we get 
\begin{equation*}
|ABC|\geq |A_{0}|+|A_{1}BC|\geq |A_{0}|+|A_{1}|+|B|=|A|+|B|,
\end{equation*}
and we are done in this case as well.

\medskip 

\noindent $(ii)$ Assume now $XC\neq X$ for every nonempty subset $X\subset AB.$

By Theorem \ref{thOl}, there exists a nonempty subset $S\subset AB$ and a
subgroup $H$ such that 
\begin{equation}
\left| S\right| \geq \left| A\right| +\left| B\right| -\left| H\right| 
\label{ineqO}
\end{equation}
and $HS=S$ or $SH=S.$ Let $c_{0}\in C$ be such that $Sc_{0}\neq S.\;$ We
claim that 
\begin{equation}
|S\cup Sc_{0}|\geq |S|+|H|.  \label{ineq}
\end{equation}

Indeed, if $HS=S,$ then $H(S\cup Sc_{0})=S\cup Sc_{0},$ and hence $S\cup
Sc_{0}$ is a disjoint union of left cosets of $H.$ This implies (\ref{ineq}) in this case.

On the other hand, if $SH=S$, pick $s_{0}\in S$ such that $s_{0}c_{0}\notin S
$, and set $H^{\prime }=H\cup \{c_{0}\}$. In the set $SH^{\prime }$, the
element $s_{0}c_{0}$ appears exactly once, since $SH\subset S.$ Applying
Theorem \ref{th_kem}, it follows that 
\begin{equation}
\left| SH^{\prime }\right| \geq \left| S\right| +\left| H^{\prime }\right|
-1=\left| S\right| +\left| H\right| .  \label{ineqST}
\end{equation}
This implies (\ref{ineq}) again, since $SH^{\prime }=S\cup Sc_{0}$.

\smallskip 

Now, as $\{1,c_{0}\}\subset C$, we have $S\cup Sc_{0}\subset SC$, and
therefore 
\begin{equation*}
\left| ABC\right| \geq |SC|\geq \left| S\cup Sc_{0}\right| .
\end{equation*}
Combined with (\ref{ineqO}) and (\ref{ineq}), this gives 
\begin{equation*}
\left| ABC\right| \geq \left| S\right| +\left| H\right| \geq \left| A\right|
+\left| B\right| .
\end{equation*}
\end{proof}

In the abelian case, the previous theorem is much shorter to prove, using
Kneser's theorem, and remains true without the assumption $B\subset C$.

\begin{theorem}
Let $A,B,C$ be nonempty finite subsets of an abelian group $G$. Suppose $1\in C.$ Then 
\begin{equation*}
ABC=AB\;\text{ or }\;\left| ABC\right| \geq \left| A\right| +\left| B\right|
.
\end{equation*}
\end{theorem}

\begin{proof}
By Kneser's theorem, we obtain 
\begin{equation*}
\left| AB\right| \geq \left| A\right| +\left| B\right| -\left| H\right| ,
\end{equation*}
where $H$ is the stabilizer of $AB$ in $G$.\ We have $HAB=AB$ and hence $HABC=ABC$. 
In particular, both $AB$ and $ABC$ are disjoint unions of cosets
of $H$. We have $AB\subset ABC$ since $1\in C.\;$ Therefore, if $ABC\neq AB$, 
it follows from the above that 
\begin{equation*}
\left| ABC\right| \geq \left| AB\right| +\left| H\right| \geq \left|
A\right| +\left| B\right| .
\end{equation*}
\end{proof}

\section{The linear setting\label{sec_linear setting}}

From now on, and for the remainder of this paper, $K$ is a commutative field
and $L$ a (possibly skew) field extension containing $K$ in its center.\
Given any subset $S\subset L,$ we write $\langle S\rangle $ for the $K$-subspace 
of $L$ generated by $S.$ For subsets $S_{1},S_{2}$ of $L$, we consider the product set 
\begin{equation*}
S_{1}S_{2}=\{s_{1}s_{2}\mid s_{1}\in S_{1},s_{2}\in S_{2}\}
\end{equation*}
and the $K$-subspace $\langle S_{1}S_{2}\rangle $ of $L$
spanned by $S_{1}S_{2}$. Note, for later use, the equality
\begin{equation*}
\langle \, \langle S_{1}\rangle \langle S_{2}\rangle \, \rangle =\langle
S_{1}S_{2}\rangle .
\end{equation*}
If $S_1 = \{x\}$, we simply write $xS_2$ instead of $\{x\}S_2$.
When $S_{1}=A$, $S_{2}=B$ are finite-dimensional $K$-subspaces of $L,$ it is easy to verify that $\langle
AB\rangle $ is finite-dimensional, with 
\begin{equation*}
\dim _{K}\langle AB\rangle \leq (\dim _{K}A)(\dim _{K}B).
\end{equation*}
Also, for any nonzero $x\in L$, the sets $xA$ and $Ax$
are $K$-subspaces of $L$, with $\dim _{K}(xA)=\dim _{K}(Ax)=\dim _{K}A.$

\medskip

\noindent \textbf{Notation.} For any subset $X \subset L$, we denote $X_{*}
= X \setminus \{0\}$ and 
\begin{equation*}
{X}_{*}^{-1} = \{x^{-1} \mid x \in X_{*} \},
\end{equation*}
the set of inverses of the nonzero elements of $X$.

\medskip

Note that ${X}_{*}^{-1}\cup \{0\}$ is not a $K$-subspace of $L$ in general, even if $X$ is.

\bigskip

In subsequent sections we establish linear versions of the addition theorems
recalled above. In particular, we obtain lower bounds on
$\dim _{K}\langle AB \rangle$ in terms of $\dim _{K}A$, $\dim _{K}B$.
As for the groups setting, our main tool will be a linear version of the
Kemperman transform.

\section{A linear Kemperman transform}

Let $(A,B)$ be a pair of
finite-dimensional $K$-subspaces of $L$. Let $x\in L\backslash \{0\}$.\ We define
the transformed pairs $(A^{\prime },B^{\prime })$ and $(A^{\prime \prime},B^{\prime \prime })$ 
with respect to $x$ as follows:
\begin{equation*}
(A^{\prime}  =  A+Ax, \, B^{\prime}  =  B\cap x^{-1}B)
\end{equation*}
and
\begin{equation*}
(A^{\prime\prime}  =  A \cap Ax^{-1},  \, B^{\prime\prime}  =  B+ xB),
\end{equation*}
where $+$ denotes the usual sum of vector subspaces.
Since $A+Ax = \langle A \cup Ax \rangle$, we may view them as linear analogues of the classical Kemperman transforms. (Compare with Section~\ref{tool}.) They satisfy the analogous properties below:

\begin{enumerate}
\item  $A^{\prime },B^{\prime },A^{\prime \prime },B^{\prime \prime }$ are $K
$-vector subspaces of $L,$

\item  $\langle A^{\prime }B^{\prime }\rangle \subset \langle AB\rangle $ \
and \ $\langle A^{\prime \prime }B^{\prime \prime }\rangle \subset \langle
AB\rangle ,$

\item  $B^{\prime }\neq \{0\}$ $\Longleftrightarrow $ $x\in B{B}_{\ast }^{-1}$ \
and \ $A^{\prime \prime }\neq \{0\}$ $\Longleftrightarrow $ $x\in {A}_{\ast
}^{-1}A,$

\item  $\dim _{K}A^{\prime }>\dim _{K}A$ $\Longleftrightarrow $ $Ax\neq A$ \
and \ \hfill \newline
$\dim _{K}B^{\prime \prime }>\dim _{K}B$ $\Longleftrightarrow $ $xB\neq B.$
\end{enumerate}

\begin{proposition}\label{one step}
\label{prop_ols} With the same notation as above, set $D={A}_{\ast }^{-1}A\cap B{B}_{\ast }^{-1}$. Suppose
that either $AD\not\subset A$ or $DB\not\subset B$. Then there exists a pair 
$(A_{1},B_{1})$ of $K$-subspaces of $L$ satisfying:

\begin{enumerate}
\item[(1)]  $A_{1}\neq \{0\}$ and $B_{1}\neq \{0\},$

\item[(2)]  $\langle A_{1}B_{1} \rangle \subset \langle AB \rangle,$

\item[(3)]  either $\dim_{K}A_{1}+\dim_{K}B_{1} = \dim_{K}A+\dim_{K}B$
and $\dim_{K}A_{1}>\dim_{K}A$, or $\dim_{K}A_{1}+\dim_{K}B_{1}>\dim_{K}A+\dim_{K}B.$
\end{enumerate}
\end{proposition}

\noindent Point (3) above looks somewhat technical, but it has a very
interesting and useful meaning. It says that, in $\mathbb{N}^2$ ordered  
lexicographically, one has
$$
\big(\dim _{K}A_{1}+\dim _{K}B_{1},\, \dim _{K}A_{1}\big) > \big(\dim _{K}A+\dim _{K}B,\, \dim _{K}A\big).
$$
This is used, for instance, in the proof of Corollary~\ref{cor2} below. Here again, note that the set $D={A}_{*}^{-1}A\cap B{B}_{*}^{-1}$ needs not be a $K$-subspace of $L$. 
\bigskip

\begin{proof}
Since $AD\not\subset A$ or $DB\not\subset B,$ there is an element $d\neq 0$
in $D$ such that $Ad\not\subset A$ or $dB\not\subset B$.\ Denote $p,q$ the dimensions
of the following quotient spaces: 
$$
p = \dim_K(A+Ad)/A, \; \; q = \dim_K(B+dB)/B.
$$
We have $\max (p,q)\geq 1$, and we need to distinguish two cases.

\smallskip
\noindent
$\bullet{}$ Assume first $p\geq q.$ In this case, we make the linear Kemperman transform
$$
A_{1}=A+Ad, \;\; B_{1}=B\cap d^{-1}B.
$$
We have $A \subsetneq A_1$ since $p \ge 1$, and $B_{1}\neq \{0\}$ since $d \in D \subset B{B}_{*}^{-1}$.
Moreover $\langle A_{1}B_{1} \rangle \subset \langle AB \rangle$. This settles the first two requirements on $(A_1,B_1)$.
As for the third one, we claim that
$$
\dim_K A_1 = \dim_K A+p, \;\; \dim_K B_1 = \dim_K B -q.
$$
Indeed, our definitions imply $p = \dim_K A_1 - \dim_K A$, and we have
\begin{eqnarray*}
\dim_K B_1 & = & \dim_K (B \cap d^{-1}B)\\
& = & \dim_K (dB \cap B)\\
& = & \dim_K dB + \dim_K B - \dim_K (dB + B)\\
& = & \dim_K B -q.
\end{eqnarray*}
It follows that
$$
\dim_{K}A_{1}+\dim_{K}B_{1}=\dim_{K}A+\dim_{K}B+p-q,
$$
whence $\dim_{K}A_{1}+\dim_{K}B_{1} \ge \dim_{K}A+\dim_{K}B$, and $\dim _{K}A_{1}>\dim _{K}A$ since $p \ge 1$.
This yields the third requirement on $(A_1,B_1)$.

\smallskip
\noindent
$\bullet{}$ Assume now $p<q.$ Here we use the other linear Kemperman transform and set 
$$
A_{1}=A\cap Ad^{-1}, \;\; B_{1}=B+dB.
$$
This time we have $B \subsetneq B_1$ since $q \ge 1$, and $A_{1}\neq \{0\}$ since $d \in D \subset {A}_{*}^{-1}A$.
We have $\langle A_{1}B_{1} \rangle \subset \langle AB \rangle$, and a similar calculation as above yields
$$
\dim_K A_1 = \dim_K A-p, \;\; \dim_K B_1 = \dim_K B+q.
$$
It follows that
$$
\dim_{K}A_{1}+\dim_{K}B_{1}=\dim_{K}A+\dim_{K}B+q-p,
$$
implying $\dim_{K}A_{1}+\dim_{K}B_{1} > \dim_{K}A+\dim_{K}B$, as desired.
\end{proof}

\begin{corollary}
\label{cor2} Let $A,B$ be nonzero finite-dimensional $K$-subspaces of $L$.
Then there exist nonzero finite-dimensional $K$-subspaces $E,F$ of $L$
satisfying
\begin{enumerate}
\item[(1)]  $\langle EF\rangle \subset \langle AB\rangle ,$
\item[(2)]  $\dim _{K}E+\dim _{K}F\geq \dim _{K}A+\dim _{K}B$,
\item[(3)] $ED=E$ and $DF=F$, where $D={E}_{\ast }^{-1}E\cap F{F}_{\ast }^{-1}.$
\end{enumerate}
\end{corollary}

\begin{proof}
By successive applications of the previous proposition, we get a
sequence 
\begin{equation*}
(A,B)=(A_{0},B_{0}),(A_{1},B_{1}),\dots ,(A_{i},B_{i}),\dots 
\end{equation*}
of pairs of $K$-subspaces of $L$ satisfying the following properties, for $i \ge 1$:
\begin{itemize}
\item $\langle A_{i}B_{i}\rangle \subset \langle A_{i-1}B_{i-1}\rangle ,$
\item $\big(\dim _{K}A_{i}+\dim _{K}B_{i},\, \dim _{K}A_{i}\big) > \big(\dim _{K}A_{i-1}+\dim _{K}B_{i-1},\, \dim _{K}A_{i-1}\big)$
\\ in $\mathbb{N}^{2}$ ordered lexicographically.
\end{itemize}
Moreover, these spaces have bounded dimension, since
\begin{eqnarray*}
\max \{\dim _{K}A_{i},\dim _{K}B_{i}\} & \leq & \dim _{K}\langle
A_{i}B_{i}\rangle  \\ 
& \leq & \dim _{K}\langle AB\rangle \\
& \leq & (\dim _{K}A)(\dim _{K}B).
\end{eqnarray*}
It follows that the above sequence must be finite. By Proposition~\ref{one step}, there is an
index $n\geq 0$ such that the set $D=(A_{n})_{\ast }^{-1}A_{n}\cap B_{n}(B_{n})_{\ast }^{-1}$
satisfies $A_{n}D\subset A_{n}$ and $DB_{n}\subset B_{n}$. In fact
$A_{n}D = A_{n}$ and $DB_{n} = B_{n}$, since $1 \in D$. Setting $E=A_n$, $F=B_n$,
it follows from the properties of the sequence of $(A_i,B_i)$ that
$\langle EF\rangle \subset \langle AB\rangle$ and $\dim _{K}E+\dim _{K}F\geq \dim _{K}A+\dim _{K}B$.
\end{proof}

\section{Linearizing Theorem~\ref{thOl} of Olson}

\label{lin ols}

With the linearized Kemperman transform at hand, we are now ready to
establish linear versions of the addition theorems of Section~\ref{sec-kem}
and \ref{variant olson}. We start with Theorem~\ref{thOl} of Olson. The
following easy lemma will be needed in the process.

\begin{lemma}
\label{lem_util} Let $D$ be a finite-dimensional $K$-subspace of $L$. Then $D
$ is a (possibly skew) field if and only if \thinspace\ $1\in D$ and $D^{2}\subset D$.
\end{lemma}

\begin{proof}
It suffices to show that every nonzero element $d\in D$ is invertible in $D$. 
Now the map $L_{d}:D\rightarrow D$ defined by $L_{d}(x)=dx$ for all $x\in D
$, is linear and injective. Hence $L_{d}$ is bijective, and therefore there
is some $d^{\prime }\in D$ with $L(d^{\prime })=1$. Obviously $d^{\prime }$
is the inverse of $d$ in $L$ and it does live in $D$.
\end{proof}

\begin{theorem}
\label{TH-Ol8lin}Let $K$ be a commutative field and $L$ a field extension of 
$K$. Let $A,B$ be finite-dimensional $K$-vector spaces in $L$ distinct
from $\{0\}$.\ Then there exist a $K$-vector subspace $S$ of $\langle AB\rangle $
and a subfield $H$ of $L$ such that

\begin{enumerate}
\item[(1)]  $K\subset H\subset L,$

\item[(2)]  $\dim _{K}S\geq \dim _{K}A+\dim _{K}B-\dim _{K}H$,

\item[(3)]  $HS=S$ or $SH=S$.
\end{enumerate}
\end{theorem}

\begin{proof}
By Corollary~\ref{cor2}, there are subspaces $E,F$ such that 
$$\langle EF\rangle  \subset \langle AB\rangle, \; \; \dim_K E + \dim_K F \ge \dim_K A + \dim_K B, $$
and $ED=E$, $DF=F$ where $D={E}_{\ast }^{-1}E\cap F{F}_{\ast }^{-1}$. 

We start by assuming $\dim_K E\geq \dim _{K}F$ and, as in the proof of Olson's
theorem, we distinguish two cases:

\smallskip
\noindent $(i)$ $F{F}_{\ast }^{-1}\not\subset {E}_{\ast }^{-1}E.$ Then there
exist $x_{1},x_{2}\in F_\ast$ such that $x_{1}x_{2}^{-1}\notin {E}_{\ast }^{-1}E$.
Therefore $Ex_{1}\cap Ex_{2}=\{0\}$ and $Ex_{1}\oplus
Ex_{2}\subset \langle EF\rangle .\;$This gives $\dim _{K}\langle EF\rangle
\geq 2\dim _{K}E\geq \dim _{K}E+\dim _{K}F.$ Let $H=K$ and $S=\langle EF\rangle $.
Then $H$ is a subfield of $L$ stabilizing $S$,
and $\dim _{K}\langle S\rangle \geq \dim _{K}A+\dim _{K}B-1$,
as desired.

\smallskip
\noindent $(ii)$ $F{F}_{\ast }^{-1}\subset {E}_{\ast }^{-1}E.$ Then 
$D=F{F}_{\ast }^{-1}$ and $1\in D$. Moreover, we have 
\begin{equation*}
DD=DF{F}_{\ast }^{-1}=F{F}_{\ast }^{-1}=D
\end{equation*}
since $DF=F.\;$So $D^{2}\subset D$.\ Let $z$ be a nonzero element of $F$.\
We have 
\begin{equation*}
F=Fz^{-1}z\subset Dz\subset DF=F
\end{equation*}
because $Fz^{-1}\subset D$ and $DF=F$.\ This implies that $Dz=F$ and thus 
$D=Fz^{-1}$ is a finite-dimensional $K$-vector space in $L$.\ Since $1\in D$
and $D^{2}=D$, we derive from Lemma \ref{lem_util} that $D$ is a subfield of 
$L$ containing $K$.\ Moreover $\dim _{K}D=\dim _{K}F$.\ Set $S=\langle
EF\rangle $ and $H=z^{-1}Dz$.\ Since $D$ is a field, $H$ is also a field and 
$\dim _{K}H=\dim _{K}D=\dim _{K}F.\;$We have 
\begin{equation*}
SH\subset \langle EFH\rangle =\langle EFz^{-1}Dz\rangle =\langle
EDzz^{-1}Dz\rangle =\langle EDz\rangle =\langle EF\rangle =S
\end{equation*}
since $Dz=F$ and $D^{2}=D.$ It follows that $SH=S$. We have $\dim _{K}S=\dim
_{K}\langle EF\rangle \geq \dim _{K}E$ since $F\neq \{0\}$.\ Finally,
this gives 
\begin{equation*}
\dim _{K}S\geq \dim _{K}E+\dim _{K}F-\dim _{K}H
\end{equation*}
because $\dim _{K}H=\dim _{K}D=\dim _{K}F.$ 

This settles the case $\dim_K E\geq \dim _{K}F$. The case 
$\dim_K F\geq \dim _{K}E$ can be treated in a similar way.
\end{proof}

\bigskip

\noindent \textbf{Remarks.}
\begin{enumerate}

\item[$(i)$] In contrast to the linear version of Kneser's
Theorem~\ref{linear Kneser}, our linear version of Olson's Theorem does not
require any separability hypothesis.

\item[$(ii)$] Assume that $L$ is commutative and $\dim_K L =p$,
a prime number. In that case, there is no intermediate field $K \subset H
\subset L$ besides $K, L$. It follows from the above theorem that either 
$\langle AB \rangle = L$ or $\dim_K \langle AB \rangle \ge \dim_K S \ge
\dim_K A + \dim_K B -1$. This is also a consequence of Theorem~\ref{linear
Kneser}, provided $L$ is further assumed to be separable over $K$.

\item[$(iii)$] Suppose that $K$ is a finite field of cardinality 
$q$.\ Then $A^{\ast }=A\backslash \{0\}$ and $B^{\ast }=B\backslash \{0\}$
are finite subsets of the group $L^{\ast }=L\backslash \{0\}$, and by
Theorem~\ref{thOl}, there exist a subset $\mathcal{S}$$^{\ast }$ of $A^{\ast
}B^{\ast}$ and a subgroup $\mathcal{H}^{\ast }$ of $L^{\ast }$ such that 
\begin{equation*}
\left| \text{$\mathcal{S}$}^{\ast }\right| \geq \left| A^{\ast }\right|
+\left| B^{\ast }\right| -\left| \text{$\mathcal{H}$}^{\ast }\right| .
\end{equation*}
Thus, we get 
\begin{equation*}
\left| A\right| +\left| B\right| \leq \left| \text{$\mathcal{S}$}^{\ast
}\right| +\left| \text{$\mathcal{H}$}^{\ast }\right| +2.
\end{equation*}
However, with Theorem \ref{TH-Ol8lin}, we obtain 
\begin{equation*}
q^{\dim _{K}A}q^{\dim _{K}B}\leq q^{\dim _{K} S}q^{\dim
_{K}H},
\end{equation*}
and since $\left| K\right| =q$, this gives 
\begin{equation*}
\left| A\right| \left| B\right| \leq \left| S\right| \left| H\right|.
\end{equation*}
Thus, Theorem \ref{thOl} gives an upper bound for $\left| A\right| +\left|
B\right|$, whereas Theorem \ref{TH-Ol8lin} gives an upper bound for $\left|
A\right| \left| B\right| $. Note that $\mathcal{S}^{\ast }\neq S\backslash
\{0\}$ and $\mathcal{H}^{\ast }\neq H\backslash \{0\}$ in general.
\end{enumerate}

\bigskip

One first easy consequence is a linear version of Kemperman's Theorem~\ref
{torsion-free} on torsion-free groups.

\begin{theorem}
\label{linear torsion-free} Let $K$ be a commutative field and $L$ a
(possibly skew) purely transcendental extension of $K$. Let $A,B$ be nonzero
finite-dimensional $K$-subspaces of $L$. Then 
\begin{equation*}
\dim _{K}\langle AB\rangle \geq \dim _{K}A+\dim _{K}B-1.
\end{equation*}
\end{theorem}

\begin{proof}
By Theorem~\ref{TH-Ol8lin}, there is a subspace $S\subset \langle AB\rangle $
and an intermediary field $K\subset H\subset L$ such that 
\begin{equation*}
\dim _{K}\langle AB\rangle \geq \dim _{K}A+\dim _{K}B-\dim _{K}H,
\end{equation*}
and $HS=S$ or $SH=S$. As $\dim _{K}S$ is finite, it follows that $\dim
_{K}H$ must be finite as well. But $K$ is the only finite-dimensional
subfield of $L$. It follows that $H=K$ and hence $\dim _{K}H=1$.
\end{proof}

\bigskip

\noindent \textbf{Remark.} The above lower bound is sharp. Indeed, let $x\in
L\setminus K$. Then $x$ is transcendental over $K$. Fix positive integers 
$r,s$. Let $A,B$ be the subspaces of $L$ generated by $\{1,x,\ldots
,x^{r-1}\} $, $\{1,x,\ldots ,x^{s-1}\}$, respectively. Then $\dim _{K}A=r$, 
$\dim _{K}B=s$ and $\dim _{K}\langle AB\rangle =r+s-1$, since 
$\langle AB\rangle $ is the subspace spanned by the basis $\{1,x,\ldots
,x^{r+s-2}\}$.

\bigskip

We now derive, from Theorem~\ref{TH-Ol8lin} again, a linear version of the
basic Proposition~\ref{basic}.

\begin{proposition}
\label{linear basic noncomm} Let $K$ be a commutative field and $L$ a field
extension containing $K$ in its center, with $\dim _{K}L$ finite. Let $A,B$
be nonzero subspaces of $L$ satisfying $\dim _{K}A+\dim _{K}B>\dim _{K}L$.
Then $\langle AB\rangle =L$.
\end{proposition}

\smallskip Note that the hypothesis that $L$ is a field is essential. For
otherwise, a counterexample would be provided by $A=L$ and $B=$ a proper
nonzero left ideal of $L$, yielding $\langle AB \rangle = B$. \medskip

\begin{proof}
By Theorem~\ref{TH-Ol8lin}, there is a subspace $S\subset \langle AB\rangle $
and an intermediate field $K\subset H\subset L$ such that 
\begin{equation}
\dim _{K}S\geq \dim _{K}A+\dim _{K}B-\dim _{K}H,  \label{olson}
\end{equation}
and either $HS=S$ or $SH=S$. We claim that $HS=L=SH$. Indeed, fix any
nonzero element $x\in L$. It follows from (\ref{olson}) and the hypothesis 
$\dim _{K}A+\dim _{K}B>\dim _{K}L$, that 
\begin{equation*}
\dim _{K}S+\dim _{K}(xH) >\dim _{K}L.
\end{equation*}
Therefore $S\cap xH\neq \{0\}$. Hence, there are nonzero elements $s\in S$
and $h\in H$ such that $s=xh$. But then $x=sh^{-1}$ belongs to $SH$, since 
$h^{-1}\in H$. It follows that $L=HS$. The same argument, with $xH$ replaced
by $Hx$, yields $L=SH$. We conclude $L=S=\langle AB\rangle $, since 
$S\subset \langle AB\rangle $ and $S=HS$ or $SH$.
\end{proof}

When $L$ is commutative, the above result admits a much simpler proof, which
does not require Theorem~\ref{TH-Ol8lin}.

\begin{proposition}
\label{linear basic comm} Let $K\subset L$ be a commutative field extension,
with $\dim _{K}L$ finite. Let $A,B$ be nonzero subspaces of $L$ satisfying 
$\dim _{K}A+\dim _{K}B>\dim _{K}L$. Then $\langle AB\rangle =L$.
\end{proposition}

\begin{proof}
We may assume $A,B\neq L$ and proceed by induction on $\dim_K B$. If $\dim_K B=1$, 
then $\dim_K A=\dim_K L$ and $\langle AB\rangle =\langle LB\rangle =L$. 
Assume now $\dim_K B\geq 2$. Since $\langle LB\rangle =L$, there must be a
nonzero element $x\in L$ such that $xB\not\subset A$. Set 
\begin{equation*}
A^{\prime }=A+xB,\; B^{\prime }=Ax^{-1}\cap B.
\end{equation*}
We have $A^{\prime }B^{\prime }\subset \langle AB\rangle $ by construction
and the commutativity of $L$. Moreover, the subspace $B^{\prime }$ is
nonzero, since $\dim _{K}Ax^{-1}+\dim _{K}B$ = $\dim _{K}A+\dim _{K}B
$ $>\dim _{K}L$. Finally, $\dim _{K}A^{\prime }+\dim _{K}B^{\prime}
=\dim _{K}A+\dim _{K}B$, since 
\begin{eqnarray*}
\dim _{K}A^{\prime }=\dim _{K}(A+xB) &=&\dim _{K}A+\dim _{K}xB-\dim
_{K}(A\cap xB) \\
&=&\dim _{K}A+\dim _{K}B-\dim _{K}(Ax^{-1}\cap B) \\
&=&\dim _{K}A+\dim _{K}B-\dim _{K}B^{\prime }.
\end{eqnarray*}
By the induction hypothesis, we conclude $\langle A^{\prime }B^{\prime
}\rangle =L=\langle AB\rangle $.
\end{proof}

\bigskip

\noindent \textbf{Remark.} We are grateful to Joseph Oesterl\'{e} for
providing us with the following alternative proof of Proposition~\ref{linear
basic noncomm}, which only uses duality in vector spaces.

\medskip

\begin{proof}
Let $H$ be any hyperplane in $L$, and let $\varphi :L\rightarrow K$ be a
linear form with kernel $H$. It suffices to show that there exist $a\in A$, 
$b\in B$ such that $\varphi (ab)\not=0$, implying $\langle AB\rangle
\not\subset H$.

The map $\beta :L\times L\rightarrow K$ defined by $\beta (x,y)=\varphi (xy)$
for all $x,y\in L$ is a $K$-bilinear form, and it is non-degenerate: if 
$x\not=0$, then $xL$ is equal to $L$ and hence $\varphi (xL)\not=\{0\}$.
Therefore $\beta $ induces an isomorphism $\gamma :L\rightarrow L^{\ast }$,
where $L^{\ast }$ is the dual of $L$, defined by the formula $\gamma
(x)(y)=\beta (x,y)$ for all $x,y\in L$.

From the inclusion map $j:B\rightarrow L$ we deduce, by transposition, a
surjection $j^{t}:L^{\ast }\rightarrow B^{\ast }$, defined as usual by 
$j^{t}(\psi )=\psi \circ j$ for all $\psi \in L^{\ast }$. The composition 
$j^{t}\circ \gamma :L\rightarrow B^{\ast }$ is then also surjective, and
therefore 
\begin{equation*}
\dim_K \ker (j^{t}\circ \gamma )=\dim_K L-\dim_K B<\dim_K A.
\end{equation*}
Thus, there is some $a\in A$ satisfying $j^{t}(\gamma (a))\not=0$. Since the
linear form $\gamma (a)\circ j:B\rightarrow K$ does not vanish, there must
be some $b\in B$ satisfying $\gamma (a)(b)\not=0$, i.e. $\beta (a,b)=\varphi
(ab)\not=0$.
\end{proof}

\section{Linearizing Theorem~\ref{THK} of Kemperman}

We shall now establish a linear analogue of Kemperman's Theorem~\ref{THK},
according to which, for subsets $A,B$ in a group $G$, one has $|AB| \ge
|A|+|B|-1$ provided there is an element $c \in AB$ with a \textit{unique
representation} of the form $c=ab$ with $a\in A, b\in B$.

In order to properly linearize this result, we need to rephrase the above
unicity condition on $c$. First, up to translation of $A,B$, we may assume
that $1 \in A \cap B$ and that $c=1$ admits the unique representation $1 =
ab $ with $a=b=1$ as a product in $AB$. If we write 
\begin{eqnarray*}
A & = & \{1\} \sqcup \overline{A} \\
B & = & \{1\} \sqcup \overline{B}
\end{eqnarray*}
with $\overline{A}, \overline{B}$ the respective complements of $\{1\}$ in $A,B$,
then the unicity of 1 as a product in $AB$ is equivalent to the disjointness
condition
\begin{equation*}
\{1\} \cap (\overline{A} \cup \overline{B} \cup \overline{A}\; \overline{B})
= \emptyset,
\end{equation*}
which may also be written as $AB = \{1\} \sqcup (\overline{A} \cup \overline{B} 
\cup \overline{A}\; \overline{B}).$ These equivalent conditions motivate
the following formulation of our sought-for linearization.

\begin{theorem}
\label{th_AB}Let $K$ be a commutative field and $L$ a field extension of $K$. 
Let $A,B$ be finite-dimensional $K$-vector spaces in $L$ such that 
$K\subset A\cap B$. Suppose there exist subspaces $\overline{A},\overline{B} \subset L$ such that 
\begin{equation}
A=K\oplus \overline{A},\text{ }B=K\oplus \overline{B}\text{ and }K\cap 
(\overline{A}+\overline{B}+\langle \overline{A}\,\overline{B}\rangle )=\{0\}.
\label{Cond}
\end{equation}
Then 
\begin{equation*}
\dim _{K}\langle AB\rangle \geq \dim _{K}A+\dim _{K}B-1.
\end{equation*}
\end{theorem}

\begin{proof}
Observe first that, if $A=K+\overline{A}$ and $B=K+\overline{B}$, then 
$\langle AB\rangle =K+(\overline{A}+\overline{B}+\langle \overline{A}\,
\overline{B}\rangle )$. Therefore, condition~(\ref{Cond}) is equivalent to 
\begin{equation*}
A=K\oplus \overline{A},\text{ }B=K\oplus \overline{B}\text{ and }\langle
AB\rangle =K\oplus (\overline{A}+\overline{B}+\langle \overline{A}\,\overline{B}\rangle ).
\end{equation*}
Since $1\in A\cap B$, we have $A+B\subset \langle AB\rangle $. We now
distinguish two cases.

\smallskip \noindent \textbf{Case 1.} Assume $A\cap B=K$. Then $\dim
_{K}(A+B)=\dim _{K}A+\dim _{K}B-1$, and we are done since $\dim _{K}\langle
AB\rangle \geq \dim _{K}(A+B)$.

\smallskip \noindent \textbf{Case 2.} Assume now $A\cap B\not=K$, i.e. $\dim
_{K}A\cap B\geq 2$. We first claim that 
\begin{equation}
A\cap B=K\oplus (\overline{A}\cap \overline{B}).
\end{equation}
Indeed, let $x\in A\cap B$. Then there are $\lambda ,\mu \in K$ and 
$\overline{x}\in \overline{A},\overline{y}\in \overline{B}$ such that 
\begin{equation*}
x=\lambda +\overline{x}=\mu +\overline{y}.
\end{equation*}
Since $K\cap (\overline{A}+\overline{B})=\{0\}$, it follows that $\lambda
=\mu $ and $\overline{x}=\overline{y}$, so that $x\in K\oplus (\overline{A}
\cap \overline{B})$, as claimed. The reverse inclusion is immediate.

Observe that in the present case, we have $\overline{A}\cap \overline{B}\not=\{0\}$. 
We shall perform a suitable sequence of linear Kemperman
transforms on the pair $(A,B)$ and eventually reach Case 1 again,
thereby concluding the proof.

Let $0\not=d\in \overline{A}\cap \overline{B}$ be any nonzero element. We
perform a Kemperman transform on $(A,B)$ relative to $d$, and get a new pair 
$(A_{1},B_{1})$, by defining either 
\begin{eqnarray}
(i) &A_{1}=A+Ad,&B_{1}=B\cap d^{-1}B,\text{ or } \\
(ii) &A_{1}=A\cap Ad^{-1},&B_{1}=B+dB.
\end{eqnarray}
In either case, we have $\langle A_{1}B_{1}\rangle \subset \langle AB\rangle 
$, and one of $(i)$ or $(ii)$ will yield the estimate 
\begin{equation*}
\dim _{K}A_{1}+\dim _{K}B_{1}\geq \dim _{K}A+\dim _{K}B.
\end{equation*}

\noindent Case $(i)$. For any nonzero element $d\in \overline{A}\cap 
\overline{B}$, define 
\begin{eqnarray*}
A_{1} &=&A+Ad, \\
B_{1} &=&B\cap d^{-1}B.
\end{eqnarray*}
We shall show that this new pair satisfies the hypotheses of the theorem.
Indeed, define 
\begin{eqnarray*}
\overline{A_{1}} &=&\overline{A}+Ad, \\
\overline{B_{1}} &=&\overline{B}\cap d^{-1}\overline{B}.
\end{eqnarray*}
We claim that conditions (\ref{Cond}) are satisfied, i.e. 
\begin{eqnarray}
A_{1} &=&K\oplus \overline{A_{1}},  \label{cond1} \\
B_{1} &=&K\oplus \overline{B_{1}},  \label{cond2} \\
\{0\} &=&K\cap (\overline{A_{1}}+\overline{B_{1}}+\langle \overline{A_{1}}\,
\overline{B_{1}}\rangle ).  \label{cond3}
\end{eqnarray}

\noindent \textbf{For (\ref{cond1}):} We have $K+\overline{A_{1}}=K+\overline{A}+Ad=A+Ad=A_{1}$. 
The sum is direct, since $\overline{A_{1}}
\subset \overline{A}+Kd+\overline{A}d\subset \overline{A}+\overline{B}
+\langle \overline{A}\,\overline{B}\rangle $, whence $K\cap \overline{A_{1}}=\{0\}$ by (\ref{Cond}).

\smallskip 

\noindent \textbf{For (\ref{cond2}):} To start with, we have $K\cap 
\overline{B_{1}}\subset K\cap \overline{B}=\{0\}$ by (\ref{Cond}).

We next verify the inclusion $K+\overline{B_{1}}\subset B_{1}$. We have 
$K\subset B$, and $K\subset d^{-1}B$ since $d\in B$, whence $K\subset B_{1}$.
Moreover, we have $\overline{B_{1}}=\overline{B}\cap d^{-1}\overline{B}
\subset B\cap d^{-1}B=B_{1}$. This establishes the desired inclusion.

It remains to prove the reverse inclusion $B_{1}\subset K+\overline{B_{1}}$.
Let $x\in B_{1}=B\cap d^{-1}B$. Since $B=K\oplus \overline{B}$, we may write 
$x=\lambda +\overline{x}$ for some $\lambda \in K$ and $\overline{x}\in 
\overline{B}$. It remains to show that $\overline{x}\in \overline{B_{1}}=
\overline{B}\cap d^{-1}\overline{B}$, i.e. that $\overline{x}\in d^{-1}
\overline{B}$. Since $x\in d^{-1}B$, there are elements $y\in B$, $\mu \in K$
and $\overline{y}\in \overline{B}$ such that $x=d^{-1}y$ and $y=\mu +
\overline{y}$. Hence $dx=y$, and therefore 
\begin{equation*}
d\lambda +d\overline{x}=\mu +\overline{y}.
\end{equation*}
It follows that $\mu =d\lambda +d\overline{x}-\overline{y}$, whence $\mu \in
K\cap (\overline{B}+\langle \overline{A}\,\overline{B}\rangle )$. Therefore 
$\mu =0$ by (\ref{Cond}), whence $d\overline{x}=\overline{y}-d\lambda \in 
\overline{B}$. This shows that $\overline{x}\in d^{-1}\overline{B}$,
implying $\overline{x}\in \overline{B_{1}}$ and finally $x\in K+\overline{B_{1}}$, as desired.

\smallskip \noindent \textbf{For (\ref{cond3}):} By (\ref{Cond}), it
suffices to show the inclusion 
\begin{equation*}
(\overline{A_{1}}+\overline{B_{1}}+\langle \overline{A_{1}}\,\overline{B_{1}}
\rangle )\subset (\overline{A}+\overline{B}+\langle \overline{A}\,\overline{B}\rangle ).
\end{equation*}
\smallskip Considering each summand at a time in the left-hand side, we have:
\newline
\smallskip $\bullet $ $\overline{A_{1}}=\overline{A}+Ad=\overline{A}+Kd+
\overline{A}d\subset (\overline{A}+\overline{B}+\langle \overline{A}\,\overline{B}\rangle )$,\newline
\smallskip $\bullet $ $\overline{B_{1}}=\overline{B}\cap d^{-1}\overline{B}\subset \overline{B},$\newline
\smallskip $\bullet $ $\overline{A_{1}}\,\overline{B_{1}}=(\overline{A}+Ad)
(\overline{B}\cap d^{-1}\overline{B})\subset \overline{A}\,\overline{B}+A\,
\overline{B}\subset \overline{A}\,\overline{B}+\overline{B}+\overline{A}\,\overline{B}$, \newline
and we are done.

\medskip 

\noindent Case $(ii)$. For any nonzero element $d\in \overline{A}\cap 
\overline{B}$, define 
\begin{eqnarray*}
A_{1} &=&A\cap Ad^{-1}, \\
B_{1} &=&B+dB.
\end{eqnarray*}
This time, we set
\begin{eqnarray*}
\overline{A_{1}} &=&\overline{A}\cap \overline{A}d^{-1}, \\
\overline{B_{1}} &=&\overline{B}+dB.
\end{eqnarray*}
With arguments similar to those of Case $(i)$, we can prove that conditions 
(\ref{Cond}) are satisfied again, i.e. 
\begin{eqnarray*}
A_{1} &=&K\oplus \overline{A_{1}}, \\
B_{1} &=&K\oplus \overline{B_{1}}, \\
\{0\} &=&K\cap (\overline{A_{1}}+\overline{B_{1}}+\langle \overline{A_{1}}\,\overline{B_{1}}\rangle ).
\end{eqnarray*}

Now, as in the proof of Corollary \ref{cor2}, we iterate the above Kemperman transforms as long as possible.
At each step, we get new subspaces $A_{i}$, $B_{i}$ satisfying (\ref
{Cond}) and such that 
$$
\big(\dim _{K}A_{i}+\dim _{K}B_{i},\, \dim _{K}A_{i}\big) > \big(\dim _{K}A_{i-1}+\dim _{K}B_{i-1},\, \dim _{K}A_{i-1}\big)
$$
in $\mathbb{N}^{2}$ ordered lexicographically.
Since these subspaces have bounded dimension, the iteration
cannot continue indefinitely and Case 1 must eventually be reached. 
This means that there exist subspaces $E,F$ of $L$ such that 
\begin{equation}
\langle EF\rangle \subset \langle AB\rangle ,  \label{H1}
\end{equation}
\begin{equation}
\dim _{K}E+\dim _{K}F\geq \dim _{K}A+\dim _{K}B,  \label{H2}
\end{equation}
and satisfying $E\cap F=K$ together with the hypotheses of the theorem.\
By Case 1, we have
\begin{equation*}
\dim _{K}\langle EF\rangle \geq \dim _{K}E+\dim _{K}F-1.
\end{equation*}
The desired inequality, namely 
\begin{equation*}
\dim _{K}\langle AB\rangle \geq \dim _{K}A+\dim _{K}B-1,
\end{equation*}
now follows from (\ref{H1}) and (\ref{H2}).
\end{proof}

\section{Linearizing Theorem~\ref{ABC}}

\label{sec: lin ABC}

Throughout the last two sections, we shall assume that $L$ is a commutative
field extension of $K$ and that every algebraic element of $L$ is separable
over $K$. This allows us to use Theorem~\ref{linear Kneser}, the linear
version of Kneser's Theorem.

Our results below probably remain true in the more general setting of the
preceding sections, where $L$ is only assumed to be a field containing $K$
in its center. But we have no proof of this so far.

\begin{lemma}
\label{Lem_util2}Let $K$ be a commutative field, $H$ a field extension of $K$
and $L$ a field extension of $H$.\ Let $V$ be a finite-dimensional $K$-vector space 
in $L$ such that $HV=V$.\ Then there exists a finite subset $R_{V}\subset V$ such that 
\begin{equation}
V=\bigoplus_{v\in R_{V}}Hv.  \label{dec}
\end{equation}
In particular $\dim _{K}H$ is finite and divides $\dim _{K}V$.
\end{lemma}

\begin{proof}
For any $v\in V$, $Hv$ is vector subspace of $V$.\ Moreover for any 
$v,v^{\prime }$ in $V$, one has $Hv=Hv^{\prime }$ or $Hv\cap Hv^{\prime
}=\{0\}.\;$The lemma follows immediately since $\dim _{K}V$ is finite.
\end{proof}

\bigskip

\noindent \textbf{Remark.} If $HV=V,$ then $V$ can be interpreted as a
finite-dimensional left $H$-module, and (\ref{dec}) gives its decomposition
into irreducible components.

\bigskip

We now give a linear version of Theorem~\ref{ABC}. As mentioned above, the
hypotheses on $L$ are probably more restrictive than actually necessary.

\begin{theorem}
\label{linear ABC} \label{Th_L}Let $K\subset L$ be commutative fields such
that every algebraic element in $L$ is separable over $K.$ Let $A,B,C\subset
L$ be finite-dimensional $K$-subspaces of $L$ such that $A,B\neq \{0\}$ and 
$K\subset C.$ Then either 
\begin{equation*}
\langle ABC\rangle =\langle AB\rangle \;\text{ or }\;\dim _{K}\langle
ABC\rangle \geq \dim _{K}A+\dim _{K}B.
\end{equation*}
\end{theorem}

\begin{proof}
We shall apply Theorem~\ref{linear Kneser}. The stabilizer $H$ of $\langle
AB\rangle $ is a field extension of $K$, and we have 
\begin{equation}
\dim _{K}\langle AB\rangle \geq \dim _{K}A+\dim _{K}B-\dim _{K}H.
\label{ineq_kneserV}
\end{equation}
Then $H$ stabilizes $\langle ABC\rangle $, since $H\langle AB\rangle
c=\langle AB\rangle c$ for all $c\in C$. Of course $\langle AB\rangle $ is a
subspace of $\langle ABC\rangle $, since $1\in C$. Assume $\langle
ABC\rangle \neq \langle AB\rangle $. Since $H\langle ABC\rangle =\langle
ABC\rangle $, Lemma~\ref{Lem_util2} implies the existence of a nonempty
finite subset $R\subset \langle ABC\rangle $ such that $R\cap \langle
AB\rangle =\emptyset $ and 
\begin{equation*}
\langle ABC\rangle =\langle AB\rangle \oplus \bigoplus_{v\in R}Hv.
\end{equation*}
In particular, we have 
\begin{equation*}
\dim _{K}\langle ABC\rangle =\dim _{K}\langle AB\rangle +\left| R\right|
\dim _{K}H.
\end{equation*}
Since $\left| R\right| >0,$ this gives 
\begin{equation*}
\dim _{K}\langle ABC\rangle \geq \dim _{K}\langle AB\rangle +\dim
_{K}H.
\end{equation*}
Finally, by (\ref{ineq_kneserV}), we obtain the desired inequality 
\begin{equation*}
\dim _{K}\langle ABC\rangle \geq \dim _{K}A+\dim _{K}B.
\end{equation*}
\end{proof}

\begin{corollary}
\label{cor3}Let $K\subset L$ be commutative fields such that every algebraic
element in $L$ is separable over $K.$ Let $A,B$ be nonzero
finite-dimensional $K$-subspaces of $L$. Then either 
\begin{equation*}
\langle AB{B}_{\ast }^{-1}B\rangle =\langle AB\rangle \;\text{ or }\;\dim
_{K}\langle AB^{2}\rangle \geq \dim _{K}A+\dim _{K}B.
\end{equation*}
\end{corollary}

\begin{proof}
Assume $\langle AB{B}_{\ast }^{-1}B\rangle \neq \langle AB\rangle $. Then
there exists $b_{0}\in B$ such that $\langle ABb_{0}^{-1}B\rangle \neq
\langle AB\rangle $.\ Thus $\langle ABb_{0}^{-1}Bb_{0}^{-1}\rangle \neq
\langle ABb_{0}^{-1}\rangle .$ We have $1\in Bb_{0}^{-1}.$ Thus, applying
the above theorem to $A$ and $Bb_{0}^{-1}$, we get 
\begin{equation*}
\dim _{K}\langle AB^{2}\rangle =\dim _{K}\langle
ABb_{0}^{-1}Bb_{0}^{-1}\rangle \geq \dim _{K}A+\dim _{K}B.
\end{equation*}
\end{proof}

\section{Powers of subspaces}

\label{sec: powers}

As in the preceding section, we assume that $L$ is a commutative field
extension of $K$ in which every algebraic element is separable over $K$. 
If $B$ is a nonzero finite-dimensional $K$-subspace of $L$, we shall consider
the sequence of powers $B$, $\langle B^2 \rangle, \langle B^3 \rangle,
\ldots $ and analyze the evolution of the nondecreasing sequence 
$\dim_K\langle B^i \rangle, i \geq 1$. Without loss of generality,
replacing $B$ by $b^{-1}B$ for some $b \in B \setminus \{0\}$, we may and
will assume that $B$ contains 1. Under this hypothesis, the sequence 
$\langle B^i \rangle$ turns into an ascending chain 
\begin{equation*}
B \subset \langle B^2 \rangle \subset \langle B^3 \rangle \subset \ldots.
\end{equation*}
This chain may eventually stabilize at $\langle B^n \rangle$ for some $n \ge
1$, for instance if $L$ is finite-dimensional over $K$. We start by
analyzing the least such exponent $n$, if any.

\begin{proposition}
\label{chain} Let $K$ be a commutative field and $L$ a field extension of $K$
containing $K$ in its center. Let $B$ be a finite-dimensional $K$-subspace
of $L$ containing 1. Let $n\geq 1$. The following are equivalent:

\begin{enumerate}
\item[(1)]  $\langle B^{n+1}\rangle =\langle B^{n}\rangle $;

\item[(2)]  $\langle B^{2n}\rangle =\langle B^{n}\rangle $;

\item[(3)]  $\langle B^{n}\rangle $ is a field.
\end{enumerate}
\end{proposition}

\begin{proof}
First observe that, if $U,V$ are any subsets of $L$, then 
\begin{equation*}
\langle UV\rangle =\langle \langle U\rangle \langle V\rangle \rangle .
\end{equation*}
Indeed, both $K$-subspaces are generated by the subset $UV$. In particular,
we may and will freely use formulas such as $\langle B^{m}\rangle =\langle
\langle B^{i}\rangle \langle B^{m-i}\rangle \rangle $ for integers $0\leq
i\leq m$.

Assume first $\langle B^{n+1}\rangle =\langle B^{n}\rangle $. Then we claim
that $\langle B^{n+i}\rangle =\langle B^{n}\rangle $ for all $i\geq 1$.
Indeed, proceeding by induction on $i$, we have 
\begin{equation*}
\langle B^{n+i}\rangle =\langle \langle B^{n+i-1}\rangle B\rangle =\langle
B^{n+1}\rangle =\langle B^{n}\rangle .
\end{equation*}
For $i=n$, this gives $\langle B^{2n}\rangle =\langle B^{n}\rangle $. In
turn, this equality is equivalent to $\langle \langle B^{n}\rangle \langle
B^{n}\rangle \rangle =\langle B^{n}\rangle $. By Lemma \ref{lem_util}, it
follows that $\langle B^{n}\rangle $ is a field.

Conversely, if $\langle B^{n}\rangle $ is a field, then $\langle
B^{2n}\rangle =\langle \langle B^{n}\rangle \langle B^{n}\rangle \rangle
=\langle B^{n}\rangle $. As $\langle B^{n}\rangle \subset \langle
B^{n+1}\rangle \subset \langle B^{2n}\rangle $, this implies that $\langle
B^{n+1}\rangle =\langle B^{n}\rangle $.
\end{proof}

\medskip

In particular, the smallest integer $n \ge 1$, if any, such that $\langle
B^{n}\rangle =\langle B^{n+1}\rangle$ coincides with the smallest integer $n
\ge 1$, if any, such that $\langle B^{n}\rangle$ is a field.

\begin{theorem}
\label{TH_last}Let $K$ be a commutative field and $L$ a commutative field
extension of $K$.$\;$Suppose that every algebraic element in $L$ is
separable over $K.$ Let $B$ be a finite-dimensional $K$-vector space in $L$
containing 1, and let $n\geq 1$. Then either 
\begin{equation*}
\langle B^{n+1}\rangle =\langle B^{n}\rangle 
\end{equation*}
or 
\begin{equation*}
\dim _{K}\langle B^{n+1}\rangle \geq \dim _{K}\langle B^{n-1}\rangle
+\dim _{K}B.
\end{equation*}
\end{theorem}

\begin{proof}
Observe first that $\langle B^{n+1}\rangle =\langle B^{n}\rangle $ if and
only if $\langle B^{n+1}\rangle =\langle B^{n}b\rangle $ for all $b\in
B\setminus \{0\}.$ Indeed, we have $\langle B^{n}\rangle \subset \langle
B^{n}b\rangle \subset \langle B^{n+1}\rangle $, and $\langle B^{n+1}\rangle $
is the sum of all the subspaces $\langle B^{n}b\rangle $ where $b$ runs over 
$B\setminus \{0\}$.

It follows that $\langle B^{n+1}\rangle =\langle B^{n}\rangle $ is
equivalent to $\langle B^{n}b_{1}\rangle =\langle B^{n}b_{2}\rangle $ for
all $b_{1},b_{2}\in B\setminus \{0\}$, which in turn is equivalent to 
$\langle B^{n}(B{B}_{\ast }^{-1})\rangle =\langle B^{n}\rangle $.

Assume now $\langle B^{n}(B{B}_{\ast }^{-1})\rangle \neq \langle
B^{n}\rangle $. Since $L$ is a commutative field, this is equivalent to 
$\langle B^{n}({B}_{\ast }^{-1}B)\rangle \neq \langle B^{n}\rangle $.
Applying Corollary \ref{cor3} to $B^{n-1}$ and $B$, this gives 
\begin{equation*}
\dim _{K}\langle B^{n+1}\rangle \geq \dim _{K}\langle B^{n-1}\rangle
+\dim _{K}B,
\end{equation*}
as required.
\end{proof}

\bigskip

\begin{corollary}
Let $K$ be a commutative field and $L$ a finite separable commutative
extension of $K$. Let $B$ be a $K$-vector space in $L$ containing 1. Then
the smallest integer $n\geq 1$ such that $\langle B^{n}\rangle $ is a field
satisfies 
\begin{equation*}
n\leq 2\,\dim _{K}L/\dim _{K}B.
\end{equation*}
\end{corollary}

\begin{proof}
By Proposition~\ref{chain}, we have $B\subsetneqq \langle B^{2}\rangle
\subsetneqq \cdot \cdot \cdot \subsetneqq \langle B^{n}\rangle $. It follows
from Theorem~\ref{TH_last} that 
$$
\dim _{K}\langle B^{n}\rangle  \geq 
\left\{
\begin{array}{rcl}
\frac{(n+1)}{2} \dim _{K}B & \; & \text{ if }n\text{ is odd}, \\
& & \\
\frac{n}{2} \dim _{K}B & \;& \text{ if }n \text{ is even}.
\end{array}
\right.
$$
Since $\dim _{K}\langle B^{n}\rangle \leq \dim _{K}L$, this imposes the
inequality 
$$
n  \le 
\left\{
\begin{array}{lll}
2\, {\dim _{K}L}/{\dim _{K}B}-1 &\;  \; & \text{ if }n\text{ is odd,}\\
& & \\
2\, {\dim _{K}L}/{\dim _{K}B}&\;  \; & \text{ if }n\text{ is even.}\\
\end{array}
\right.
$$
\end{proof}

\bigskip

\noindent \textbf{Remark.} From the previous Corollary, we deduce 
\begin{equation*}
n\leq \left\lfloor \frac{2\dim _{K}L}{\dim _{K}B}\right\rfloor.
\end{equation*}
This upper bound is sharp.\ This can easily be seen, for example by choosing
for $B$ a supplementary space of $K$ in $L.$

\bigskip \smallskip \noindent \textbf{Acknowledgment.} The authors thank
Vincent Fleckinger, Joseph Oesterl\'{e} and Surya D. Ramana for stimulating
discussions concerning Section~\ref{lin ols} of this paper.


\begin{thebibliography}{9}
\bibitem{Xian}  \textsc{X. D.\ Hou, K.\ H. Leung and Xiang. Q,} \textit{A
generalization of an addition theorem of Kneser}, J. Number Theory 
\textbf{97} (2002), 1-9.

\bibitem{Kem}  \textsc{J. H.\ B.\ Kemperman,} \textit{On complexes in a
semigroup}, Indag.\ Math\ \textbf{18} (1956), 247-254.

\bibitem{Kn53}  \textsc{M.\ Kneser,} \textit{Absch\"{a}tzung der
asymptotischen Dichte von Summenmengen}, Math. Z. \textbf{58} (1953),
459-484.

\bibitem{Kn55} \textsc{M.\ Kneser,} \emph{Ein Satz \"uber abelsche Gruppen mit Anwendungen auf die Geometrie der Zahlen}, 
Math. Z. \textbf{61} (1955), 429-434. 

\bibitem{ols}  \textsc{J. E.\ Olson,} \textit{On the sum of two sets in a
group}, J. Number Theory \textbf{18} (1984), 110-120.
\end{thebibliography}
\end{document}